\documentclass[12pt]{article}
\usepackage{amssymb}
\usepackage{latexsym}
\usepackage{amsmath}

\newtheorem{theorem}{Theorem}[section]
\newtheorem{lemma}[theorem]{Lemma}

\newtheorem{definition}[theorem]{Definition}
\newtheorem{example}[theorem]{Example}
\newtheorem{corollary}[theorem]{Corollary}
\newtheorem{proposition}[theorem]{Proposition}
\newtheorem{remark}[theorem]{Remark}

\newcommand {\ld}{\underline{d}}
\newcommand {\ud}{\overline{d}}
\newcommand {\lbd}{\underline{B\!D}}
\newcommand {\ubd}{\overline{B\!D}}
\newcommand {\lld}{\underline{ld}}
\newcommand {\uld}{\overline{ld}}
\newcommand {\z}{\mathbb{Z}}
\newcommand {\n}{\mathbb{N}}

\newcommand{\I}{\mathcal{I}}

\title{Abstract densities and ideals of sets}

\author{Mauro Di Nasso\thanks
{Supported by PRIN 2012 ``Models and Set'', MIUR.}\\
\normalsize{\texttt{mauro.di.nasso@unipi.it}}
\and
Renling Jin\thanks
{Renling Jin would like to thank the Universit\`a di Pisa and
PRIN 2012 ``Models and Set'', MIUR, for supporting this collaborative research.}\\
\normalsize{\texttt{jinr@cofc.edu}}
}

\begin{document}

\date{}

\maketitle

\begin{abstract}
Abstract upper densities are monotone and subadditive
functions from the power set of positive integers
to the unit real interval that generalize the
upper densities used in number theory, including
the upper asymptotic density, the upper Banach density, and the
upper logarithmic density.

We answer a question posed by G. Grekos in 2013,
and prove the existence of translation invariant
abstract upper densities onto the unit interval,
whose null sets are precisely
the family of finite sets, or the family of
sequences whose series of reciprocals converge.
We also show that no such density can be atomless.
(More generally, these results also hold for a large class
of summable ideals.)


\end{abstract}

\section*{Introduction}

Several notions of densities for sets of natural numbers
are used in number theory for different purposes, including
the upper and lower asymptotic densities, the upper and lower
Banach (or uniform) density, the upper and lower logarithmic
density, the Schnirelmann density, \emph{etc.}.
The idea of an abstract notion of density that encompasses
the basic features of the known densities have been repeatedly
considered in the literature (see, \emph{e.g.}
\cite{CB53,A67,M76,FS81,ST83,G05}).

At the open problem session
of the Workshop \emph{``Densities and their applications''},
held in St. Etienne in July 2013,
the following question was asked by G. Grekos:

\begin{itemize}
\item
\emph{Is there an ``abstract density''
$\delta$ such that $\delta(A)=0$ if and only if $A$
is finite? Or such that $\delta(A)=0$ if and only if
$\sum_{a\in A}1/a<\infty$?
More generally,
given an ideal of sets
$\mathcal{I}\subseteq\mathcal{P}(\n)$,
is there an ``abstract density'' $\delta$ such that
$\delta(A)=0$ if and only if $A\in\mathcal{I}$?}
\end{itemize}

Recall that a nonempty family
$\mathcal{I}\subseteq\mathcal{P}(\n)$
is an \emph{ideal} if it closed under subsets
and under taking finite unions, and $\n\notin\I$.

To avoid trivial examples, such an abstract density
should behave ``nicely'', in the sense that it should share as many
as possible of the properties of the familiar densities
as considered in number theory.
In this paper we investigate around the above questions.

The paper is organized as follows.
In the first section, we discuss the general properties
of an abstract density; in the second section
we present our main results;
% in a conventional sense;
the third section contains the proofs.

\medskip
\section{Abstract Densities}

Let $\n$ be the set of positive integers.
Let us start by isolating the fundamental
features that an abstract notion of density on $\n$ must have.

\begin{definition}\label{densitydef}
{\rm An \emph{abstract density} on $\n$ is a function
$\delta:\mathcal{P}(\n)\to[0,1]$ defined on the family
of all subsets of $\n$, taking values in the unit real interval,
and that satisfies the following properties:

\begin{enumerate}
\item[(1)]
$\delta(\n)=1$,
\item[(2)]
$\delta(F)=0$ for every finite $F\subset\n$,
\item[(3)]
\emph{Monotonicity}:
If $A\subseteq B$ then $\delta(A)\le\delta(B)$.
\end{enumerate}}
\end{definition}

We remark that virtually all upper and lower
densities that have been considered
in number theory are examples of abstract densities in the above sense.
For example, it is easily seen that the following seven densities
for subsets of $\n$ satisfy properties (1), (2), (3).\footnote
{~Also other densities that have been considered
in the literature satisfy properties
(1), (2), (3), \emph{e.g.}
upper Buck, upper
Polya, upper analytic, and exponential densities
(see \cite{LT15} for the definitions).}

For $a\leq b$ in $\n$, we denote
$A(a,b):=|A\cap [a,b]|$ and by $A(a):=A(1,a)$.

\begin{itemize}
\item
The \emph{lower asymptotic density}:
\[\ld(A):=\liminf_{n\rightarrow\infty}\frac{A(n)}{n};\]
\item
The \emph{upper asymptotic density}:
\[\ud(A):=\limsup_{n\rightarrow\infty}\frac{A(n)}{n};\]
\item
The \emph{lower Banach density}:
\[\lbd(A):=\lim_{n\rightarrow\infty}\inf_{k\in\n}\frac{A(k+1,k+n)}{n};\]
\item
The \emph{upper Banach density}
(or simply \emph{Banach density} or \emph{uniform density}):
\[\ubd(A):=\lim_{n\rightarrow\infty}\sup_{k\in\n}\frac{A(k+1,k+n)}{n};\]
\item
The \emph{lower logarithmic density}:
\[\lld(A):=\liminf_{n\rightarrow\infty}\frac{
\sum_{a\in A\cap[1,n]}a^{-1}}{\sum_{a=1}^n a^{-1}};\]
\item
The \emph{upper logarithmic density}:
\[\uld(A):=\limsup_{n\rightarrow\infty}\frac{
\sum_{a\in A\cap[1,n]}a^{-1}}{\sum_{a=1}^na^{-1}};\]
\item
The \emph{Shnirelmann density}:\footnote
{~Schnirelmann density is included here because
it is widely used in number theory.
However, it is worth remarking that it is
a ``one of a kind'' density, that was designed
to guarantee the implication:
``$\sigma(A)=1\Rightarrow A=\n$.''}
\[\sigma(A)=\inf_{n\geq 1}\frac{A(n)}{n}.\]
\end{itemize}

Since here we are interested in abstract densities $\delta$
whose family of null sets is closed under finite unions,
also the following property should be satisfied:
$$``\delta(A)=\delta(B)=0\ \Rightarrow\ \delta(A\cup B)=0."$$
In consequence, a natural assumption is subadditivity.

\begin{definition}
{\rm An abstract density is called \emph{abstract upper density}
if it satisfies:
\begin{enumerate}
\item[(4)]
\emph{Subadditivity}:
$\delta(A\cup B)\le\delta(A)+\delta(B)$.
\end{enumerate}}
\end{definition}

The name ``upper density'' is justified by the
fact that the three upper densities itemized above
(namely, upper asymptotic density, upper Banach
density, and upper logarithmic density) are indeed
subadditive. On the contrary, the three corresponding
lower densities, as well as Shnirelmann density, are not.

The following is easily proved:

\begin{itemize}
\item
\emph{Let $\delta$ be an abstract upper density.
If the symmetric difference $A\bigtriangleup B$ is finite,
then $\delta(A)=\delta(B)$.}
\end{itemize}

The monotonicity and subadditivity properties
are independent of each other.

\begin{proposition}
Let $\delta:\mathcal{P}(\n)\to[0,1]$. Then
\begin{itemize}
\item[(i)]
Properties (1), (2), (3) do not imply (4);
\item[(ii)]
Properties (1), (2), (4) do not imply (3).
\end{itemize}
\end{proposition}

\smallskip
\noindent
\textsc{Proof.}
$(i)$.
The lower asymptotic density is an example
of an abstract density that is not subadditive.
(\emph{E.g.}, let $A=\bigcup_{n\in\n}(a_{2n},a_{2n+1}]$
where $\{a_n\}$ is any increasing sequence with
$\lim_{n\to\infty}\frac{a_{n+1}}{a_n}=0$;
then it is easily verified that $\underline{d}(A)=\underline{d}(A^c)=0$,
while $\underline{d}(A\cup A^c)=1$.)

$(ii)$.
Let $E\subseteq\n$ be the set of all even numbers, let
$O\subseteq\n$ be the set of all odd numbers, let
$E_1$ be the set of multiples of $4$, and let
$E_2=E\setminus E_1$
be the set of numbers that are congruent
to $2$ modulo 4. Define the
function $\delta:\mathcal{P}(\n)\to[0,1]$
as follows:
$$\delta(A)\ =\
\begin{cases}
0 & \mbox{if } A \text{ is finite,}
\\
1 & \mbox{if } A\cap O\ \mbox{is infinite,}
\\
1 & \mbox{if } A\cap O\ \mbox{is finite, and exactly one of }
A\cap E_1\ \mbox{and } A\cap E_2\ \mbox{is infinite,}
\\
\frac{1}{2} & \mbox{if } A\cap O\ \mbox{is finite, and both }
A\cap E_1\ \mbox{and } A\cap E_2\ \mbox{are infinite.}
\end{cases}$$

It is readily seen that $\delta$ satisfies (1) and (2).
Moreover, $\delta$ is subadditive. Indeed,
if both $A$ and $B$ are infinite, then
$\delta(A\cup B)\leq 1=
\frac{1}{2}+\frac{1}{2}\leq\delta(A)+\delta(B)$;
and if at least one of the two sets is finite, say $B$,
then $\delta(A\cup B)=\delta(A)\leqslant\delta(A)+\delta(B)$.
However, $\delta$ is not monotonic; for example,
$E_1\subseteq E$ but $\delta(E_1)=1>\delta(E)=1/2$.
\qquad $\Box$

\bigskip
For $A\subseteq\n$ and $m\in\n$, denote by
$A+m:=\{a+m\mid a\in A\}$ the \emph{rightward translation}
of $A$ by $m$; and by
$A-m:=\{k\in\n\mid k+m\in A\}$ the \emph{leftward translation}
of $A$ by $m$.
A natural property to be considered for densities
is the following:

\begin{itemize}
\item[(5)]
\emph{Translation invariance}:
$\delta(A+m)=\delta(A-m)=\delta(A)$ for every $A\subseteq\n$ and
for every $m\in\n$.
\end{itemize}

It is easily seen that
an abstract upper density that is rightward translation invariant,
is also leftward translation invariant. Indeed, in this case,
$\delta(A-m)=\delta((A-m)+m)=\delta(A)$, because
the symmetric difference $((A-m)+m)\bigtriangleup A$
is finite.

With the only exception of Shnirelmann density,
we remark that all other densities itemized at the
beginning of this section are translation invariant
(see, \emph{e.g.}, \cite{LT15}).

There is a trivial example of an abstract density that fulfills
all properties (1)--(5).

\begin{example}\label{zero-one}
Let $\delta:\mathcal{P}(\n)\to[0,1]$ be the function
defined by setting $\delta(A)=0$ if $A$ is finite,
and $\delta(A)=1$ if $A$ is infinite.
Then $\delta$ is a translation invariant upper abstract density.
\end{example}

The above example calls for additional
properties that make abstract densities
``non-trivial''. Two such properties
that have been repeatedly considered in
the literature are the following:

\begin{enumerate}
\item[(6)]
\emph{Richness}:
For every $r\in [0,1]$ there exists $A\subseteq\n$ with $\delta(A)=r$;
\item[(7)]
\emph{Atomless-ness}:
For every $A\subseteq\n$ with $\delta(A)>0$, there exists
$B\subset A$ such that $0<\delta(B)<\delta(A)$.
\end{enumerate}

Notice that even by assuming (1) -- (5), richness and
atomless-ness are independent of each other.
The following example is motivated by \cite{O77}.

\begin{example}\label{atomlessnotrich}
For every $A\subseteq\n$, let
$$\delta(A)\ :=\
\begin{cases}
\frac{1}{2}(1+\ud(A)) & \mbox{ if }\,\ud(A)>0,
\\
0 & \mbox{ if }\,\ud(A)=0.
\end{cases}$$
Then $\delta$ is a translation invariant
abstract upper density that is atomless but not rich, because
$\text{range}(\delta)=\{0\}\cup\left(\frac{1}{2},1\right]$.
\end{example}

The notion of asymptotic density also makes sense
when relativized to any infinite set.

\begin{definition}
{\rm Let $X\subseteq\n$ be an infinite set. For every $A\subseteq\n$,
the \emph{upper asymptotic density of $A$ relative to} $X$ is
defined as
\[\ud_X(A):=\limsup_{n\rightarrow\infty}\frac{(A\cap X)(n)}{X(n)}.\]}
\end{definition}

Notice that relative upper densities are
abstract upper densities that are rich;
however, in general, they are not translation invariant.

\begin{example}
Let $\langle b_n\mid n\in\n\rangle$
be a ``rapidly growing" sequence of natural numbers,
in the sense that $\lim_{n\to\infty}\frac{b_{n+1}-b_n}{b_n}=\infty$.
Let $B:=\bigcup_{n=1}^{\infty}[b_n,b_{n+1}]$,
and for $A\subseteq\n$ let
$$\delta_B(A):=\begin{cases}
1 & \mbox{ if }\ud_B(A)>0,
\\
0 & \mbox{ if }\ud_B(A)=0
\end{cases}
\quad\text{and}\quad
\delta(A):=\max\{\delta_B(A),\ud_{\n\setminus B}(A)\}.$$
Then $\delta$ is
a translation invariant abstract upper density that is rich
but not atomless. Indeed, by the choice of $B$, both
$\ud_B$ and $\ud_{\n\setminus B}$ are translation invariant.
Moreover, $\delta$ is rich because $\ud_{\n\setminus B}$ is rich.
Finally, $B$ is an ``atom", that is,
$\delta(B)=1$ and for every $B'\subseteq B$ one has either
$\delta(B')=0$ or $\delta(B')=\delta(B)$.
\end{example}

\smallskip
A stronger property that directly implies both
richness and atomless-ness is the following intermediate
value property:

\begin{itemize}
\item[(8)]
\emph{Darboux property}:\footnote
{~The name ``Darboux property'' has been used in the literature
because it resembles the intermediate value property of derivatives,
as stated in Darboux's theorem in real analysis
(see the remarks in \cite[\S 2]{LT17}).}
For every $A\subseteq\n$ and for every $0\le r\le \delta(A)$, there exists
a subset $B\subseteq A$ such that $\delta(B)=r$.
\end{itemize}

\bigskip
\section{The results}\label{main}

Recall that an \emph{ideal} (over $\n$) is a nonempty
family $\mathcal{I}\subseteq\mathcal{P}(\n)$ with
$\n\not\in\mathcal{I}$, which is
closed under taking finite unions and subsets; that is
$A,B\in\mathcal{I}\Rightarrow A\cup B\in\mathcal{I}$, and
$B\subseteq A\in\mathcal{I}\Rightarrow B\in\mathcal{I}$.

We say that an ideal $\I\subseteq\mathcal{P}(\n)$
is \emph{translation invariant}
if $A\in\I\Leftrightarrow A+k\in\I$ for every $k\in\z$.

The first easy example of a translation invariant
ideal is given by the family of finite sets:
$$\textsf{fin}:=\{A\subseteq\n\mid A\mbox{ is finite}\,\}.$$

Another relevant example of translation invariant ideal
is given by the
family of those sequences whose series of reciprocals converge:
$$\textsf{rcp}:=
\Big\{A\subseteq\n\ \Big|\,\sum_{a\in A}1/a<\infty\Big\}.$$

Recall the following general notion (\emph{e.g.}, see \cite{F00}).
For every non-increasing function $f:\n\to[0,\infty)$
with $\sum_{n=1}^\infty f(n)=\infty$,
the family
$$\I_f\ =\ \left\{A\subseteq\n\,\Big|\,\sum_{n\in A}f(n)<\infty\right\}$$
is the \emph{summable ideal} determined by $f$.
It is easily verified that such a family $\I_f$ is
indeed a translation invariant ideal that includes
all finite sets.

Notice that both $\textsf{fin}=\I_f$ and
$\textsf{rcp}=\I_g$ are summable ideals,
where $f$ is the constant function with value $1$,
and $g(n)=1/n$ is the ``reciprocal'' function, respectively.

Abstract densities and ideals are closely related notions;
indeed, as one can easily verify, the family
of \emph{zero sets}
$$\mathcal{Z}_{\delta}\ :=\ \{A\subseteq\n\mid\delta(A)=0\}$$
of any (translation invariant) abstract upper density $\delta$ is a
(translation invariant) ideal over $\n$.

The following question was posed by G. Grekos at the open
problem session of the Workshop
\emph{``Densities and their applications''},
held in St. Etienne in July 2013:

\bigskip
\noindent
\textbf{Question. (G. Grekos)}\ \
Given an ideal $\mathcal{I}$ on $\n$,
for example $\mathcal{I}=\textsf{fin}$ or
$\mathcal{I}=\textsf{rcp}$,
does there exist a ``nice density'' $\delta$ such that
$\mathcal{Z}_\delta=\mathcal{I}$ ?

\bigskip
In this section we present the main results
that we obtained in this paper
to address the above question.
Proofs of theorems will be given in the next section.

\bigskip
Two sets $A,B\subseteq\n$
are called $\I$-\emph{almost disjoint} ($\I$-AD for short)
if $A\cap B\in\I$; and are called
$\I$-\emph{translation almost disjoint}
($\I$-TAD for short) if for every $s,t\in\z$,
the translates $A+s,B+t$ are $\I$-almost disjoint.
When $\I$ is translation invariant, the latter condition
is equivalent to having $(A+k)\cap B\in\I$ for every $k\in\z$.

The above notions are extended to families of sets
in a natural way: A family of infinite sets
$\mathcal{A}\subseteq\mathcal{P}(\n)$
is $\I$-AD (or $\I$-TAD) if every pair $A,B$ of distinct elements of
$\mathcal{A}$ are $\I$-AD (or $\I$-TAD, respectively).
In order to make definitions more meaningful,
it is also assumed that such families $\mathcal{A}$ do not
contain any member of $\I$.

Following the usual terminology, we will simply say
``almost disjoint'' and ``translation almost disjoint'' to mean
$\textsf{fin}$-AD and $\textsf{fin}$-TAD, respectively.

\begin{theorem}\label{one}
Let $\mathcal{I}\subseteq\mathcal{P}(\n)$
be a translation invariant ideal that includes all finite sets, and
assume that there exists an infinite $\mathcal{I}$-TAD family.
Then there exists an abstract upper density
$\delta$ that is translation invariant and rich,
and such that $\mathcal{Z}_\delta=\mathcal{I}$.
\end{theorem}

\begin{theorem}\label{two}
Let $\I_f$ be the summable ideal determined by $f$.
Then there exists
an $\I_f$-TAD family $\mathcal{A}$ of the cardinality
of the continuum. In consequence,
there exists an abstract upper density
$\delta$ that is translation invariant and rich,
and such that $\mathcal{Z}_\delta=\I_f$.
\end{theorem}

\begin{corollary}\label{corollaryoftwo}
There exist abstract upper densities $\delta_1$ and $\delta_2$
such that $\z_{\delta_1}=\textsf{fin}$ and
$\z_{\delta_2}=\textsf{rcp}$, respectively,
that are translation invariant and rich.
\end{corollary}

On the negative side, one cannot ask
for the abstract densities of Corollary \ref{corollaryoftwo}
to be also atomless, and even more so, to satisfy Darboux property.

We say that $A$ is \emph{$\I$-almost included} in $B$,
and we write $A\subseteq_\I B$, if $A\setminus B\in\I$.
An ideal $\I$ has the \emph{diagonal intersection property}
(DIP for short) if for all sequences $\langle B_n\mid n\in\n\rangle$
where $B_n\notin\I$ and $B_{n+1}\subseteq_\I B_n$ for every $n$,
there exists a set $A\notin\I$ such
that $A\subseteq_\I B_n$ for all $n$.

\begin{proposition}\label{DIP}
If $\mathcal{I}_f$ is a summable ideal determined by a
 non-increasing function $f$, then $\mathcal{I}_f$
satisfies the DIP.
\end{proposition}

\textsc{Proof}.
Suppose that $\langle B_n\mid n\in\n\rangle$
is a sequence of sets of natural numbers where
$B_n\notin\mathcal{I}_f$
and $B_{n+1}\setminus B_n\in\mathcal{I}_f$ for all $n$.
Let $a_1=\min B_1$ and, proceeding by induction,
assume that elements $a_1<a_2<\ldots<a_{n(s)}$ have been
found such that $a_1,a_2,\ldots,a_{n(t)}\in\bigcap_{i=1}^t B_i$
and $\sum_{i=1}^{n(t)}f(a_i)\geq t$ for each $t\leq s$.
Since $\bigcap_{i=1}^{s+1}B_i\notin\mathcal{I}_f$,
there exist elements $a_{n(s)+1}<\ldots< a_{n(s+1)}$
with $a_{n(s)+1}>a_{n(s)}$ such that
$a_{n(s)+1},\ldots, a_{n(s+1)}\in \bigcap_{i=1}^{s+1}B_i$
and $\sum_{i=n(s)+1}^{n(s+1)}f(a_i)\geq 1$.
Let $A=\{a_i\mid i\in\n\}$.
Clearly, $A\not\in\mathcal{I}_f$ because
$\sum_{i=1}^{n(s)}f(a_i)\geq s$ for all $s$ implies
that $\sum_{a\in A}f(a)=\infty$;
besides, for every $s$ we have that
$A\setminus B_s\subseteq\{a_i\mid i=1,\ldots, n(s)\}\in
\textsf{fin}\subset\mathcal{I}_f$.
\qquad $\Box$

\begin{corollary}
Both the ideal $\textsf{fin}$ and the ideal $\textsf{rcp}$
satisfy the DIP.
\end{corollary}

\begin{theorem}\label{three}
If an abstract upper density $\delta$ is atomless
then the ideal $\mathcal{Z}_\delta$ of its zero sets
does not satisfy the DIP.
\end{theorem}

\begin{corollary}\label{corollaryofthree}
If $\delta$ is an atomless abstract upper density, then
$\mathcal{Z}_{\delta}\not=\textsf{fin}$
and $\mathcal{Z}_{\delta}\not=\textsf{rcp}$.
\end{corollary}

\section{The proofs}\label{proofs}

\textsc{Proof of Theorem \ref{one}}.\
Let $\mathcal{A}_0$ be an infinite $\I$-TAD family.
By a direct application of \emph{Zorn's Lemma},
we can pick a maximal $\I$-TAD
family $\mathcal{A}\supseteq\mathcal{A}_0$.
%Clearly such an $\mathcal{A}$ is infinite.
Enumerate its elements $\mathcal{A}=\{A_\alpha\mid \alpha<\mu\}$,
where $\mu=|\mathcal{A}|$ is an infinite cardinal.
Notice that for every $m\in\n$ and for every $k\in\z$,
the translation invariance of $\I$ guarantees that
$(A+k)\cap B\in\I$ if and only if $(A+k+m)\cap(B+m)\in\I$;
in consequence, $A$ and $B$ are $\I$-TAD if and only if
$A$ and $B+m$ are $\I$-TAD.
Notice also that, by maximality,
for every $B\notin\I$
there must be $A_{\alpha}\in\mathcal{A}$
such that $A_{\alpha}$ and $B$ are not $\I$-TAD.

We are now ready to
construct the desired abstract upper density $\delta$.
For $n\in\n$ and $B\subseteq\n$, we set
$\delta_n(B):=0$ if $A_n$ and $B$ are $\I$-TAD;
otherwise we set
$$\delta_n(B):=\frac{1}{n+1}+\frac{1}{n(n+1)}\cdot
\sup_{k\in\z}\left(\ud_{A_n+k}(B)\right)\in
\left[\frac{1}{n+1},\frac{1}{n}\right].$$

In case $\mu>\aleph_0$, for infinite ordinals
$\alpha<\mu$
we set $\delta_\alpha(B)=0$ if $A_\alpha$ and $B$
are $\I$-TAD; otherwise we set $\delta_\alpha(B)=1$.
Finally, we define $\delta:\mathcal{P}(\n)\to[0,1]$ by letting:
$$\delta(B):=\sup_{\alpha<\mu}\delta_\alpha(B).$$

Let us now verify that $\delta$ satisfies the required properties.
Notice first that $\delta(B)=0$ for every $B\in\I$.
Indeed, for every $\alpha$, it directly follows from the
the definition of $\delta_\alpha$ that $\delta_\alpha(B)=0$
whenever $B\in\I$.

All pairs $A_\alpha$ and $\n$ are
not $\I$-TAD, since
$(A_\alpha+k)\cap\n=A_\alpha+k\notin\I$ for every $k\in\z$.
(Here we used the facts that $\mathcal{A}\cap\I=\emptyset$
and that $\I$ is translation invariant.)
Now, trivially $\ud_{A_1}(\n)=1$, and so
$\delta_1(\n)=1/2+1/2\cdot 1=1$,
and hence $\delta(\n)=1$.

If $B\subseteq B'$, then it is readily verified that
$\delta_\alpha(B)\le\delta_\alpha(B')$ for every $\alpha$.
By passing those inequalities to the limit superior
as $\alpha<\mu$, we obtain the monotonicity
property $\delta(B)\le\delta(B')$.

Next, we show that for every $\alpha$
and for every $B,C\subseteq\n$,
one has the inequality
$\delta_\alpha(B\cup C)\le\delta_\alpha(B)+\delta_\alpha(C)$.
Clearly, this will prove the subadditivity of $\delta$, because
\begin{eqnarray}
\nonumber
\delta(B\cup C) & = &
\sup_{\alpha<\mu}\delta_\alpha(B\cup C)\ \le\
\sup_{\alpha<\mu}(\delta_\alpha(B)+\delta_\alpha(C))
\\
\nonumber
{} & \le &
\sup_{\alpha<\mu}\delta_\alpha(B)+\sup_{\alpha<\mu}\delta_\alpha(C)\ =\
\delta(B)+\delta(C).
\end{eqnarray}
If $n\in\n$, the desired inequality for $\delta_n$
follows directly from the definition. Indeed,
$$\delta_n(B\cup C)\ \le\ \frac{1}{n}\ \le\
\frac{1}{n+1}+\frac{1}{n+1}\ \le\ \delta_n(B)+\delta_n(C).$$
When $\alpha$ is infinite,
$\delta_\alpha$ only assumes
the values 0 and 1. So, if by contradiction
$\delta(B\cup C)>\delta_\alpha(B)+\delta_\alpha(C)$, then it
must be $\delta_\alpha(B\cup C)=1$ and
$\delta_\alpha(B)=\delta_\alpha(C)=0$.
But this is impossible because if both pairs $A_\alpha,B$ and
$A_\alpha,C$ are $\I$-TAD,
then also the pair $A_\alpha,B\cup C$ is $\I$-TAD.

The density $\delta$ is translation invariant because
for every $\alpha$, for every $B\subseteq\n$, and
for every $m\in\n$, we have $\delta_\alpha(B+m)=\delta_\alpha(B)$.
Indeed, if $n\in\n$ then for every $k\in\z$
one has $\ud_{A_n+k}(B)=\ud_{A_n+k+m}(B+m)$, so
$\sup_{k\in\z}\left(\ud_{A_n+k}(B)\right)=
\sup_{h\in\z}\left(\ud_{A_n+h}(B+m)\right)$,
and hence $\delta_n(B)=\delta_n(B+m)$.
If $\alpha$ is infinite, $\delta_\alpha(B)=\delta_\alpha(B+m)$
because
$A_\alpha$ and $B$ are $\I$-TAD if and only if $A_\alpha$
and $B+m$ are $\I$-TAD. Indeed, by the translation
invariance of $\I$, for every $k\in\z$
one has that $(A+k)\cap B\in\I$ if and only
if $[(A+k)\cap B]+m=(A+k+m)\cap(B+m)\in\I$.

Let us now turn to the richness property.
Clearly, $\delta(B)=0$ whenever $B\in\I$.
Given $r\in(0,1]$, pick $n_0\in\n$ and $\lambda\in[0,1]$
such that $r=\frac{1}{n_0+1}+\frac{1}{n_0(n_0+1)}\cdot\lambda$.
Since $\ud_{A_{n_0}}$ is rich,
there exists an infinite $B\subseteq A_{n_0}$ such that
$\ud_{A_{n_0}}(B)=\lambda$.
Trivially $(A_{n_0}\cap B)(m)=B(m)$ for every $m\in\n$
and so, for every $k\in\z$ one has
$$\frac{((A_{n_0}+k)\cap B)(m)}{(A_{n_0}+k)(m)}\ \le\
\frac{B(m)}{A_{n_0}(m)-|k|}\ =\
\frac{(A_{n_0}\cap B)(m)}{A_{n_0}(m)}\cdot
\frac{A_{n_0}(m)}{A_{n_0}(m)-|k|}.$$
By passing to the limit superiors as $n$ goes to infinity,
we get $\ud_{A_{n_0}+k}(B)\le\ud_{A_{n_0}}(B)$,
and hence
$\sup_{k\in\mathbb{Z}}\ud_{A_{n_0}+k}(B)=\ud_{A_{n_0}}(B)=\lambda$.
For $\alpha\ne n_0$,
the sets $A_\alpha$ and $A_{n_0}$ are $\I$-TAD, and hence also
the sets $A_\alpha$ and $B$ are $\I$-TAD.
Then $\delta_\alpha(B)=0$ for $\alpha\ne n_0$, and
so $\delta(B)=\delta_{n_0}(B)=r$, as desired.

We have noticed already that if $B\in\I$ then $\delta(B)=0$,
so let us assume that $B\notin\I$.
By the maximality of the $\I$-TAD family $\{A_\alpha\mid\alpha<\mu\}$,
there exists $\alpha$ such that $A_\alpha$ and $B$ are not $\I$-TAD.
If such an $\alpha=n$ is finite, then
$\delta(B)\ge\delta_n(B)\ge\frac{1}{n+1}>0$;
and if such an $\alpha$ is infinite, then $\delta(B)=\delta_\alpha(B)=1$.
This shows that $\mathcal{Z}_\delta=\I$.
\qquad $\Box$

\bigskip
\noindent
\textsc{Proof of Theorem \ref{two}}.\
We prove a lemma first.

\begin{lemma}\label{functiong}
Let $\I_f$ be the summable ideal determined by $f$.
Then there exist increasing functions
$g,h:\n\to\n$ and a sequence $n_0<n_1<n_2<\cdots$ such that
\begin{enumerate}
\item The sequences $\{g(n)-g(n-1)\}_{n\in\n}$ and
$\{h(n)-h(n-1)\}_{n\in\n}$ are non-decreasing,
\item $\lim_{n\rightarrow\infty}(g(n)-g(n-1))=\infty$,
\item $h(n_m)-h(n_m-1)\geq 2^m$ for every positive integer $m$,
\item $\sum_{n=1}^\infty f(g(h(n)))=\infty$.
\end{enumerate}
\end{lemma}

\noindent
\emph{Proof of Lemma.}
We first define by induction
a function $g(n)$ and a
sequence $l_0<l_1<l_2<\cdots$ so that
$g(l_m)-g(l_m-1)\geq 2^m$. Set $l_0=0$ and $g(0)=0$.
Suppose we have found $g(0)<g(1)<\cdots<g(l_m)$
such that
\begin{itemize}
\item $\{g(i)-g(i-1)\}_{i\leq l_m}$ is non-decreasing,
\item $g(l_m)-g(l_m-1)\geq 2^m$, and
\item $\sum_{i=0}^{l_m}f(g(i))\geq m$.
\end{itemize}
We now define
$l_{m+1}$ and $g(i)$ for $i=l_m+1,l_m+2,\ldots,l_{m+1}$.

Let $d=g(l_m)-g(l_m-1)$. Since $\sum_{i>g(l_m)}f(i)=\infty$,
there exists an $i_0<2d$ such that
$\sum_{k=2}^{\infty}f(g(l_m)+i_0+2dk)=\infty$ because $[g(l_m),\infty)$ is the
disjoint union
of $2d$ arithmetic sequences of common difference $2d$.
Let $l_{m+1}$ be sufficiently large such that
$\sum_{k=2}^{l_{m+1}}f(g(l_m)+i_0+2dk)\geq 1$. We set
\begin{itemize}
\item $g(l_m+1)=g(l_m)+i_0$ and $g(l_m+2)=g(l_m)+i_0+2d$
if $i_0\geq d$ or
\item $g(l_m+1)=g(l_m)+d$ and $g(l_m+2)=g(l_m)+i_0+2d$ if $i_0<d$.
\end{itemize}
We also set $g(l_m+j+1)=g(l_m+j)+2d$ for $j=2,3,\ldots,l_{m+1}$.
It is easy to check that
\begin{itemize}
\item $\{g(i)-g(i-1)\}_{i\leq l_{m+1}}$ is non-decreasing,
\item $g(l_{m+1})-g(l_{m+1}-1)\geq 2d\geq 2^{m+1}$, and
\item $\sum_{i=0}^{l_{m+1}}f(g(i))\geq m+1$.
\end{itemize}
Notice that the purpose of choosing the value of
$g(l_m+1)$ to be $g(l_m)+i_0$ or to be $g(l_m)+d$ is to guarantee that 
the sequence $\{g(i)-g(i-1)\}_{i\leq l_{m+1}}$ be non-decreasing.

Clearly, $\lim_{n\rightarrow\infty}(g(n)-g(n-1))=\infty$
because $g(l_m)-g(l_m-1)\geq 2^m$.

Now we define the function $h(n)$ and the sequence $\{n_m\}$
exactly as above, by replacing $f(n)$ with $f(g(n))$.
This completes the proof of Lemma \ref{functiong}\qquad $\Box$

\bigskip

Let us go back to the proof of Theorem \ref{two}.

Assume that the function $g,h$ and the sequence $\{n_m\}$ are
constructed as in Lemma \ref{functiong}.
Let $2^{\n}:=\{\sigma\mid\sigma:\n\rightarrow\{0,1\}\}$
and $2^{[1,n]}:=\{\sigma\mid\sigma:[1,n]\rightarrow\{0,1\}\}$. Notice that
$|2^{[1,n]}|=2^n$.
We construct a collection of sets
$\mathcal{A}=\{A_{\sigma}\subseteq G\mid\sigma\in 2^{\n}\}$ where
$G=\mbox{range}(g)$
such that for any distinct $\sigma,\tau\in 2^{\n}$,
\begin{itemize}
\item $|A_{\sigma}\cap [g(h(n-1)),g(h(n))-1]|=1$ for all $n\in\n$,
\item $A_{\sigma}\cap A_{\tau}\subseteq [0,g(h(n_{m}-1))]$ where
$m=\min\{n\mid\sigma(n)\not=\tau(n)\}$.
\end{itemize}

The theorem follows from the construction of $\mathcal{A}$:
Clearly, $\mathcal{A}$ is an AD of the cardinality of the continuum.
Since $n$-th element of $A_{\sigma}$ is between $g(h(n-1))$
and $g(h(n))-1$, we have that
$\sum_{a\in A_{\sigma}}f(a)\geqslant\sum_{n=2}^{\infty}f(g(h(n)))=\infty$,
which implies $A_{\sigma}\not\in \I_f$. If $x\in G\cap (G+k)$,
then $x,x-k\in G$, and so $G\cap (G+k)$ must
be finite because $\lim_{n\rightarrow\infty}(g(n)-g(n-1))=\infty$.
Therefore, $G$ and $G+k$ are almost disjoint for any non-zero $k\in\mathbb{Z}$.
As a consequence, $A\cap (B+k)$ is a finite set for any
distinct pair $A,B\in\mathcal{A}$ and any $k\in\mathbb{Z}$.
Thus $\mathcal{A}$ is a TAD.

We now construct $\mathcal{A}$. For any positive integer
$m$ and $\sigma\in 2^{\n}$,
let $\sigma\!\upharpoonright\! m$ represent the restriction of the function
$\sigma$ on $\{1,2,\ldots,m\}$, let 
$A_{\sigma\upharpoonright m}:=A_{\sigma}\cap
[0,g(h(n_m))-1]$, and let $\mathcal{A}\!\upharpoonright\! m:=
\{A_{\sigma\upharpoonright m}\mid \sigma\in 2^{\n}\}$.
We construct $\mathcal{A}$ by defining
$\mathcal{A}\!\upharpoonright\! m$ inductively on $m$.

Let $A_{\sigma\upharpoonright 1}=\{0\}$. 
Now suppose we have obtained
$\mathcal{A}\!\upharpoonright\! m$. Since
\[\{h(n)-h(n-1)\}_{n\in\n}\] is non-decreasing, we have that
$h(n_m+i)-h(n_m-1+i)\geq h(n_m)-h(n_m-1)\geq 2^m$ 
for $0\leq i< n_{m+1}-n_m$. Thus there are
at least $2^m$-many distinct elements in each of the sets
  \[\{G\cap [g(h(n_m-1+i)),g(h(n_m+i))-1]\mid 0\leq i<n_{m+1}-n_m\}.\]
  Let $\{a_{i,\gamma}\mid \gamma\in 2^{[1,m]}\}$ be an enumeration of 
  a set of cardinality $2^m$ in \[G\cap [g(h(n_m-1+i)),g(h(n_m+i))-1]\]
  for $i=1,2,\ldots,n_{m+1}-n_m-1$
  and $A'_{\sigma\upharpoonright m}=A_{\sigma\upharpoonright m}\cup
  \{a_{i,\sigma\upharpoonright m}\mid i=1,2,\ldots,n_{m+1}-n_m-1\}$.
  It is easy to see that
  \begin{itemize}
  \item $A_{\sigma\upharpoonright m}
  =A'_{\sigma\upharpoonright m}\cap [0,g(h(n_m))-1]$;
  \item $A'_{\sigma\upharpoonright m}\cap A'_{\tau\upharpoonright m}
  =A_{\sigma\upharpoonright m}\cap A_{\tau\upharpoonright m}$
  for any $\sigma\!\upharpoonright\! m\not=\tau\!\upharpoonright\! m$;
  \item $A'_{\sigma\upharpoonright m}\cap [g(h(n_m-1+i),g(h(n_m+i))-1]
  =\{a_{i,\sigma\upharpoonright m}\}$ for $0<i<n_{m+1}-n_m$.
  \end{itemize}
Since $h(n_{m+1})-h(n_{m+1}-1)\geq 2^{m+1}$, we can
label $2^{m+1}$ distinct elements in $G\cap [g(h(n_{m+1}-1)),g(h(n_{m+1}))-1]$ by
$\{b_{\gamma}\mid \gamma\in 2^{[1,m+1]}\}$. Now let
$A_{\sigma\upharpoonright (m+1)}=A'_{\sigma\upharpoonright m}\cup\{b_{\gamma}\}$
where $\gamma=\sigma\!\upharpoonright\!(m+1)$. Notice that
each set $A'_{\sigma\upharpoonright m}$ has two different extensions
in $\mathcal{A}\!\upharpoonright\!(m+1)$ depending on the value of 
$\sigma(m+1)$.

This completes the construction of $\mathcal{A}\!\upharpoonright\!(m+1)$.
Now let $A_{\sigma}=\bigcup_{m=1}^{\infty}A_{\sigma\upharpoonright m}$.
The set $\mathcal{A}=\{A_{\sigma}\mid
\sigma\in 2^{\n}\}$ is the desired AD family.\qquad $\Box$

\bigskip

If $\I=\textsf{fin}$ or $\I=\textsf{rcp}$, the constructions of
$g$ and $h$ in Lemma \ref{functiong} can be simplified.
In fact $g$ and $h$ can be the same function.

It is readily verified that the sequence $g(n)=h(n)=n^2$
satisfies the required properties for
the summable ideal $\textsf{fin}=\I_f$,
where $f$ is the constant function with value $1$.

For the summable ideal $\textsf{rcp}=\I_f$ where $f(n)=1/n$ let
$\varphi(n)=\lfloor n\cdot\log\log n\rfloor$, where
$\lfloor\,\cdot\,\rfloor$ denotes the integer part.
Since $\varphi(\varphi(n))=o(n\log n)$ and
$\sum_{n=2}^\infty \frac{1}{n\log n}=\infty$,
we have that $\sum_{n=3}^\infty \frac{1}{\varphi(\varphi(n))}=\infty$.
Moreover, it is easily checked that
$\{\varphi(n+1)-\varphi(n)\}_{n\ge 16}$
is a non-decreasing and unbounded sequence of natural numbers.
So, the previous theorem applies
by taking $g(n)=h(n)=\varphi(n+15)$.
\qquad $\Box$

\begin{remark}
We do not know whether $g$ and $h$ in Lemma \ref{functiong}
can be the same function for a generic summable ideal $I_f$.
\end{remark}

\bigskip
\noindent
\textsc{Proof of Theorem \ref{three}}.\
For an upper abstract density $\delta$,
let us write $A\subseteq_\delta B$ if
$A$ is included in $B$ up to a set of zero density, that is,
when $\delta(A\setminus B)=0$.
Clearly, $A\subseteq_\delta B$ implies that $\delta(A)\le\delta(B)$.
The desired result directly follows from the following

\begin{lemma}
Let $\delta$ be an upper abstract density.
If $\delta$ is atomless then there exists a decreasing sequence
$\langle B_n\mid n\in\n\rangle$ of sets of positive
density such that the following holds:
For every $A\subseteq\n$, if $A\subseteq_\delta B_n$ for
all $n$, then $\delta(A)=0$.
\end{lemma}

\noindent
\emph{Proof of Lemma.}
For each $X\subseteq\n$ with $\delta(X)>0$, define
$$\gamma(X)\ =\
\inf\{\delta(B)\mid B\subseteq X\,\mbox{ and }\,\delta(B)>0\}.$$
Since $\delta$ is atomless, $\gamma(X)<\delta(X)$.
Moreover, $\gamma$ is non-increasing, that is,
if $Y\subseteq X$ and $\delta(Y)>0$, then
$\gamma(X)\leq\gamma(Y)$.
Now let $B_1=\n$.
At the inductive step,
let $\eta_n=\delta(B_n)-\gamma(B_n)>0$, and let
$\varepsilon_n=\frac{1}{2}(\gamma(B_n)+\delta(B_n))>\gamma(B_n)$.
Then we can pick a subset $B_{n+1}\subset B_n$
such that $0<\delta(B_{n+1})\le\varepsilon_n$.
Since $\gamma(B_{n+1})\ge\gamma(B_n)$, we have
$$\eta_{n+1}\ =\ \delta(B_{n+1})-\gamma(B_{n+1})\ \le\
\varepsilon_n-\gamma(B_n)\ =\ \frac{\eta_n}{2}.$$
In consequence, $\lim_{n\to\infty}\eta_n=0$, and so
$\lim_{n\to\infty}\gamma(B_n)=\lim_{n\to\infty}\delta(B_n)$.
If $A\subseteq_\delta B_n$ for all $n$
then $\gamma(B_n)\le\delta(A)\le\delta(B_n)$ for all $n$, and hence
$\delta(A)=\lim_{n\to\infty}\gamma(B_n)=\lim_{n\to\infty}\delta(B_n)$.
If by contradiction $\delta(A)>0$, we could pick $A'\subset A$
such that $0<\delta(A')<\delta(A)$, and
we would have $\delta(A')<\gamma(B_n)$
for all but finitely many $n$. This is not possible because
$0<\delta(A'\cap B_n)=\delta(A')<\gamma(B_n )$, against
the definition of $\gamma(B_n)$.
\qquad $\Box$

\end{document}